\documentclass[a4paper,11pt]{amsart}

\usepackage{graphicx}
\usepackage{mathptmx}
\usepackage{amsmath}
\usepackage{amssymb}
\usepackage{enumitem}
\usepackage{xcolor}

\newmuskip\pFqmuskip

\newcommand*\pFq[6][8]{%
  \begingroup 
  \pFqmuskip=#1mu\relax
  \mathcode`=\string"8000
  \begingroup\lccode`\~=`\,
  \lowercase{\endgroup\let~}\pFqcomma
  F^{#2}_{#3}{\left(\genfrac..{0pt}{}{#4}{#5}\bigg|#6\right)}%
  \endgroup
}
\newcommand{\pFqcomma}{\mskip\pFqmuskip}

\newtheorem{theorem}{Theorem}[section]

\newtheorem{corollary}[theorem]{Corollary}

\begin{document}

\title[]{Probabilistic degenerate Stirling numbers of the first kind and their applications}

\author{Taekyun Kim}
\address{Department of Mathematics, Kwangwoon University, Seoul 139-701, Republic of Korea}
\email{tkkim@kw.ac.kr}
\author{Dae San Kim}
\address{Department of Mathematics, Sogang University, Seoul 121-742, Republic of Korea}
\email{dskim@sogang.ac.kr}


\begin{abstract}
Let $Y$ be a random variable whose degenerate moment generating functions exist in some neighborhoods of the origin. The aim of this paper is to study the probabilistic degenerate Stirling numbers of the first kind associated with $Y$ which are constructed from the degenerate cumulant generating function of $Y$. They are a degenerate version of the probabilistic Stirling numbers of the first kind associated with $Y$, which were recently introduced by Adell-B\'enyi. We investigate some properties, related identities, recurrence relations and explicit expressions for those numbers. In addition, we apply our results to the special cases of normal and gamma random variables.
\end{abstract}

\maketitle

\markboth{\centerline{\scriptsize Probabilistic degenerate Stirling numbers of the first kind and their applications}}
{\centerline{\scriptsize T. Kim and D. S. Kim}}

\section{Introduction} 
Let $Y$ be a random variable whose degenerate moment generating functions of $Y$ exist in some neighborhoods of the origin (see \eqref{10-1}).
In this paper, we study the probabilistic degenerate Stirling numbers of the first kind associated with $Y$, $S_{1,\lambda}^{Y}(n,k)$, which are constructed from the degenerate cumulant generating function of $Y$ (see \eqref{19}, \eqref{21}). They are a degenerate version of the probabilistic Stirling numbers of the first kind associated with $Y$, $S_{1}^{Y}(n,k)$ (see \eqref {16}), which were recently introduced by Adell-B\'enyi (see [2]). The aim of this paper is to investigate some properties, related identities, recurrence relations and explicit expressions for $S_{1,\lambda}^{Y}(n,k)$. Furthermore, we apply our results to normal random variable with parameters $(\mu,\sigma^{2})=(1,1)$ and gamma random variable with parameters $(\alpha, \beta)=(0,1)$. Here we note that various degenerate versions, certain $\lambda$-analogues and probabilistic extensions of many special numbers and polynomials have been explored and some remarkable results have been obtained (see [2,3,10,16-18,20-22,26 and the references therein]). \par
The outline of this paper is as follows. In Section 1, we recall the degenerate exponentials (see \eqref{1}), the degenerate logarithms (see \eqref{6}), the Stirling numbers of the first kind, $S_{1}(n,k)$ (see \eqref{2}) and the degenerate Stirling numbers of first kind, $S_{1,\lambda}(n,k)$ (see \eqref{4}). We remind the reader of the gamma random variable $Y \sim \Gamma(\alpha,\beta)$ with parameters $\alpha, \beta$, the normal random variable $Y \sim N(\mu,\sigma^{2})$ with parameters $\mu, \sigma^{2}$ and the numbers $\kappa_{n}(Y)$ in terms of the cumulant generating function (see \eqref{10}). Let $Y$ be as before. Then we recall the probabilistic degenerate Stirling numbers of the second kind associated with $Y$, ${n \brace k}_{Y,\lambda}$ (see \eqref{11}, \eqref{12}), and the probabilistic Stirling numbers of the first kind associated with $Y$, $S_{1}^{Y}(n,k)$ (see \eqref{16}), which were introduced by Adell-B\'enyi in [2]. We remind the reader of the probabilistic degenerate Bernoulli polynomials associated with $Y$, $\beta_{n,\lambda}^{Y}(x)$ (see \eqref{13}) and the probabilistic degenerate Euler polynomials associated with $Y$, $\mathcal{E}_{n,\lambda}^{Y}(x)$ (see \eqref{15}). We recall the Lah numbers $L(n,k)$ and the degenerate unsigned Stirling numbers of the first kind. Section 2 contains the main results of this paper. We first introduce $\kappa_{n,\lambda}(Y)$, as a degenerate version of $\kappa_{n}(Y)$, in terms of the degenerate cumulant generating function of $Y$ (see \eqref{19}). In Theorem 2.1, we express $\kappa_{n,\lambda}(Y)$ as a finite sum involving ${n \brace k}_{Y,\lambda}$. Then we define the probabilistic degenerate Stirling numbers of the first kind associated with $Y$, $S_{1,\lambda}^{Y}(n,k)$, as a degenerate version of $S_{1}^{Y}(n,k)$. In Theorem 2.2, we find two different expressions for $\big(E[e_{\lambda}^{Y}(t)]\big)^{x}$, one involving ${n \brace k}_{Y,\lambda}$ and the other involving $S_{1,\lambda}^{Y}(n,k)$. In Theorem 2.4, we express ${n \brace l}_{Y,\lambda}$ in terms of ${k \brace l}_{-\lambda}$, and $S_{1,\lambda}^{Y}(n,l)$ in terms of $S_{1,-\lambda}(k,l)$. In Theorem 2.5, we express ${n \brace j}_{Y, \lambda}$, as a finite sum involving $S_{1,\lambda}^{Y}(n,k)$, $L(k,l)$ and ${l \brace j}_{\lambda}$. We represent $E[(Y)_{n,\lambda}]$ as a finite sum involving $S_{1,\lambda}^{Y}(n,k)$ and $L(k,l)$ in Theorem 2.6. We find recurrence relations for $S_{1,\lambda}^{Y}(n,k)$ in Theorem 2.7 and for ${n \brace k}_{Y,\lambda}$ in Theorem 2.9, which give respectively the values for $S_{1,\lambda}^{Y}(k,k)$ and ${k \brace k}_{Y,\lambda}$. In Theorem 2.10, we represent $\beta_{n,\lambda}^{Y}(x)$ in terms of $\beta_{l,\lambda}^{Y}=\beta_{l,\lambda}^{Y}(0)$ and $S_{1,\lambda}^{Y}(j,k)$. Similar result is obtained for $\mathcal{E}_{n,\lambda}^{Y}(x)$ in Theorem 2.11. In Theorem 2.12 we express $\frac{1}{2}\big(\sum_{l=0}^{n}\binom{n}{l}\mathcal{E}_{l,\lambda}^{Y}(x)E\big[(Y)_{n-l,\lambda}\big]+\mathcal{E}_{n,\lambda}^{Y}(x)\big)$ as a finite sum involving $S_{1,\lambda}^{Y}(m,k)$ and ${n \brace m}_{Y, \lambda}$. In Section 3, we apply our results to normal and gamma random variables. Assume that $Y \sim N(0,1)$ for Theorems 3.1-3.4. Then we show that $E[e_{\lambda}^{Y}(t)]=e^{\frac{1}{2}(\log e_{\lambda}(t))^{2}}$, from which we express $E[(Y)_{n,\lambda}]$ in terms of $S_{1}(m,2k)$ and ${n \brace m}_{\lambda}$. In Theorem 3.2, we represent ${n \brace k}_{Y, \lambda}$ as a finite sum involving $S_{1}(m,2j)$ and ${n \brace m}_{\lambda}$. We express $\mathcal{E}_{n,\lambda}^{Y}(x)$ in terms of $E_{k}(x)$, $S_{1}(j,2k)$ and ${n \brace j}_{\lambda}$ in Theorem 3.3. In Theorem 3.4, we represent $\frac{1}{2}\big(\sum_{l=0}^{n}\binom{n}{l}\mathcal{E}_{l,\lambda}^{Y}(x)E\big[(Y)_{n-l,\lambda}\big]+\mathcal{E}_{n,\lambda}^{Y}(x)\big)$ as a polynomial in powers of $x$ with coefficients involving $S_{1}(j,2k)$ and ${n \brace j}_{\lambda}$. Assume that $Y \sim \Gamma(1,1)$ for Theorems 3.5 and 3.6. In Theorem 3.5, we represent $S_{1,\lambda}^{Y}(n,k)$ as a finite sum involving $S_{1}(n,l)$ and $S_{1,\lambda}(l,k)$, which reduces to $S_{1}(n,k)$ as $\lambda$ tends to 0. Finally, in Theorem 3.6 we express $\frac{1}{m+1}\Big(\beta_{m+1,\lambda}^{Y}(n+1)-\beta_{m+1,\lambda}^{Y}\Big)$ as a finite sum involving $S_{1}(m,l)$. In the rest of this section, we recall the facts that are needed throughout this paper.

\vspace{0.1in}

For any nonzero $\lambda\in\mathbb{R}$, the degenerate exponentials are defined by 
\begin{equation}
e_{\lambda}^{x}(t)=\sum_{n=0}^{\infty}(x)_{n,\lambda}\frac{t^{n}}{n!},\quad e_{\lambda}^{1}(t)=e_{\lambda}(t),\quad (\mathrm{see}\ [3,9,15-21]),\label{1}	
\end{equation}
where $(x)_{0,\lambda}=1,\ (x)_{n,\lambda}=x(x-\lambda)(x-2\lambda)\cdots(x-(n-1)\lambda),\ (n\ge 1)$. 
Note that $\displaystyle\lim_{\lambda\rightarrow 0}e_{\lambda}^{x}(t)=e^{xt}$. For $n\ge 0$, the Stirling numbers of the first kind are defined by 
\begin{equation}
(x)_{n}=\sum_{k=0}^{n}S_{1}(n,k)x^{k},\quad (\mathrm{see}\ [1-26]),\label{2}	
\end{equation}
where $(x)_{0}=1,\ (x)_{n}=x(x-1)\cdots(x-n+1),\ (n\ge 1)$. \par 
The Stirling numbers of the second kind are given by 
\begin{equation}
x^{n}=\sum_{k=0}^{n}{n\brace k}(x)_{k},\quad (n\ge 0),\quad (\mathrm{see}\ [11,18,23]). \label{3}
\end{equation}
Recently, the degenerate Stirling numbers of the first kind are defined by Kim-Kim as
\begin{equation}
(x)_{n}=\sum_{k=0}^{n}S_{1,\lambda}(n,k)(x)_{k,\lambda},\quad (n\ge 0),\quad (\mathrm{see}\ [16]).\label{4}
\end{equation}
As the inverse relation of \eqref{4}, the degenerate Stirling numbers of the second kind are defined by 
\begin{equation}
(x)_{n,\lambda}=\sum_{k=0}^{n}{n\brace k}_{\lambda}(x)_{k},\quad (n\ge 0),\quad (\mathrm{see}\ [16]). \label{5}
\end{equation}
Note that $\displaystyle\lim_{\lambda\rightarrow 0}{n\brace k}_{\lambda}={n\brace k}$ and $\lim_{\lambda\rightarrow 0}S_{1,\lambda}(n,k)=S_{1}(n,k)$. \\
The degenerate logarithm $\log_{\lambda}t$ is defined as the compositional inverse of $e_{\lambda}(t)$, so that 
\begin{equation}
\log_{\lambda}(1+t)=\frac{1}{\lambda}\big((1+t)^{\lambda}-1\big)=\sum_{n=1}^{\infty}\lambda^{n-1}(1)_{n,1/\lambda}\frac{t^{n}}{n!},\quad (\mathrm{see}\ [15-22]). \label{6}
\end{equation}
Note that $\lim_{\lambda\rightarrow 0}\log_{\lambda}(1+t)=\log(1+t)$. We recall here that the generating function of the degenerate Stirling numbers of the first kind is given by
\begin{equation}
\frac{1}{k!}\big(\log_{\lambda}(1+t)\big)^{k}=\sum_{n=k}^{\infty}S_{1,\lambda}(n,k)\frac{t^{n}}{n!}. \label{7}
\end{equation}
The continuous random variable $Y$ whose density function is given by 
\begin{equation}
f(y)=\left\{\begin{array}{ccc}
	\beta e^{-\beta y}\frac{(\beta y)^{\alpha-1}}{\Gamma(\alpha)}, & \textrm{if $y\ge 0$}, \\
	0, & \textrm{if $y<0$},
\end{array}\right.\label{8}
\end{equation}
for some $\alpha,\beta>0$, is said to be the gamma random variable with parameters $\alpha,\beta$, and denoted by $Y\sim\Gamma(\alpha,\beta)$, (see [18,22,24]). \par 
A continuous random variable $Y$ is said to be the  normal random variable with parameters $\mu,\sigma^{2}$, denoted by $Y\sim N(\mu,\sigma^{2})$, if its probability density function is given by 
\begin{equation}
f(y)=\frac{1}{\sqrt{2\pi}\sigma}e^{-(x-\mu)^{2}/2\sigma^{2}},\quad \mathrm{for}\quad y\in(-\infty,\infty),\quad (\mathrm{see}\ [18,21,24]).\label{9}
\end{equation}
The cumulant generating function of $Y$ is defined by 
\begin{equation}
K_{Y}(t)=\log\Big(E\big[e^{Yt}\big]\Big)=\sum_{n=0}^{\infty}\kappa_{n}(Y)\frac{t^{n}}{n!},\quad (\mathrm{see}\ [2,24]). \label{10}	
\end{equation}
Thus, by \eqref{10}, we get 
\begin{displaymath}
\kappa_{1}(Y)=E[Y],\quad \kappa_{2}(Y)=\mathrm{Var}(Y),\quad \kappa_{3}(Y)=E\big[(Y-E(Y))^{3}\big].
\end{displaymath}
Throughout this paper, we assume that $Y$ is a random variable such that the degenerate moment generating functions of $Y$, 
\begin{equation}
E\big[e_{\lambda}^{Y}(t)\big]=\sum_{n=0}^{\infty}E\big[(Y)_{n,\lambda}\big]\frac{t^{n}}{n!},\quad (|t|<r=r(\lambda)), \label{10-1}
\end{equation}
exist for some $r=r(\lambda)>0$. Let $(Y_{j})_{j\ge 1}$ be a sequence of mutually independent copies of the random variable $Y$, and let $S_{k}=Y_{1}+Y_{2}+\cdots+Y_{k}, (k\ge 1)$, with $S_{0}=0$. Then the probabilistic degenerate Stirling numbers of the second kind associated with $Y$ are defined by 
\begin{equation}
\frac{1}{k!}\sum_{j=0}^{k}\binom{k}{j}(-1)^{k-j}E\big[(S_{j})_{n,\lambda}\big]={n\brace k}_{Y,\lambda},\quad (n\ge k\ge 0),\quad (\mathrm{see}\ [18,21,22]).\label{11}
\end{equation}
Alternatively, they are given by
\begin{equation}
\frac{1}{k!}\Big(E\big[e_{\lambda}^{Y}(t)\big]-1\Big)^{k}=\sum_{n=k}^{\infty}{n \brace k}_{Y,\lambda}\frac{t^{n}}{n!},\quad (k\ge 0). \label{12}	
\end{equation}
In [22], the probabilistic degenerate Bernoulli polynomials associated with $Y$ are defined by 
\begin{equation}
\frac{t}{E\big[e_{\lambda}^{Y}(t)\big]-1}\Big(E\big[e_{\lambda}^{Y}(t)\big]\Big)^{x}=\sum_{n=0}^{\infty}\beta_{n,\lambda}^{Y}(x)\frac{t^{n}}{n!}.\label{13}	
\end{equation}
Note that $\displaystyle\lim_{\lambda\rightarrow 0}\beta_{n,\lambda}^{Y}(x)=B_{n}^{Y}(x)$, where $B_{n}^{Y}(x)$ are the probabilistic Bernoulli polynomials associated with $Y$ given by 
\begin{equation}
\frac{t}{E\big[e^{Yt}\big]-1}\Big(E\big[e^{Yt}\big]\Big)^{x}=\sum_{n=0}^{\infty}B_{n}^{Y}(x)\frac{t^{n}}{n!},\quad (\mathrm{see}\ [17,22]). \label{14}
\end{equation}
The probabilistic degenerate Euler polynomials associated with $Y$ are given by 
\begin{equation}
\frac{2}{E\big[e_{\lambda}^{Y}(t)\big]+1}\Big(E\big[e_{\lambda}^{Y}(t)\big]\Big)^{x}=\sum_{n=0}^{\infty}\mathcal{E}_{n,\lambda}^{Y}(x)\frac{t^{n}}{n!},\quad (\mathrm{see}\ [22]). \label{15}
\end{equation}
In [2], the probabilistic Stirling numbers of the first kind associated with $Y$, $S_{1}^{Y}(n,k)$, are defined by 
\begin{equation}
\frac{1}{k!}\Big(\log E\big[e^{Yt}\big]\Big)^{k}=\sum_{n=k}^{\infty}(-1)^{n-k}S_{1}^{Y}(n,k)\frac{t^{n}}{n!},\quad (k\ge 0). \label{16}
\end{equation}
It is well known that the Lah numbers are defined by 
\begin{equation}
\langle x\rangle_{n}=\sum_{k=0}^{n}L(n,k)(x)_{k},\quad (n\ge 0),\quad (\mathrm{see}\ [11,12,23]), \label{17}
\end{equation}
where $\langle x\rangle_{0}=1,\ \langle x\rangle_{n}=x(x+1)\cdots(x+n-1),\ (n\ge 1)$. \par 
The degenerate unsigned Stirling numbers of the first kind are defined by 
\begin{equation}
	\langle x\rangle_{n}=\sum_{k=0}^{n}{n \brack k}_{\lambda}\langle x\rangle_{k,\lambda}, (n\ge k\ge 0),\quad (\mathrm{see}\ [15,16,17,20]). \label{18}
\end{equation}
Note that ${n\brack k}_{\lambda}=(-1)^{n-k}S_{1,\lambda}(n,k),\ (n\ge k\ge 0)$. 

\section{Probabilistic degenerate Stirling numbers of the first kind} 
Now, we consider the {\it{degenerate cumulant generating function of $Y$}} which is defined by 
\begin{equation}
\log_{-\lambda}\Big(E\big[e_{\lambda}^{Y}(t)\big]\Big)=\sum_{n=0}^{\infty}\kappa_{n,\lambda}(Y)\frac{t^{n}}{n!}.\label{19}	
\end{equation}
Note that $\displaystyle\lim_{\lambda\rightarrow 0}\kappa_{n,\lambda}(Y)=\kappa_{n}(Y),\ (n\ge 0)$. From \eqref{19}, we note that 
\begin{align}
\sum_{n=0}^{\infty}\kappa_{n,\lambda}(Y)\frac{t^{n}}{n!}&=\log_{-\lambda}\Big(E\big[e_{\lambda}^{Y}(t)\big]-1+1\Big)\label{20}\\
&=\sum_{k=1}^{\infty}\frac{(-1)^{k-1}\lambda^{k-1}}{k!}(1)_{k,-1/\lambda}\Big(E\big[e_{\lambda}^{Y}(t)\big]-1\Big)^{k}\nonumber \\
&=\sum_{k=1}^{\infty}(-1)^{k-1}\lambda^{k-1}\langle 1\rangle_{k,1/\lambda}\sum_{n=k}^{\infty}{n \brace k}_{Y,\lambda}\frac{t^{n}}{n!}\nonumber \\
&=\sum_{n=1}^{\infty}\sum_{k=1}^{n}(-1)^{k-1}\lambda^{k-1}\langle 1\rangle_{k,1/\lambda}{n\brace k}_{Y,\lambda}\frac{t^{n}}{n!},\nonumber
\end{align}
where $\langle x\rangle_{0,\lambda}=1,\ \langle x\rangle_{n,\lambda}=x(x+\lambda)(x+2\lambda)\cdots(x+(n-1)\lambda),\ (n\ge 1)$. Therefore, by comparing the coefficients on both sides of \eqref{20}, we obtain the following theorem. 
\begin{theorem}
Let $n$ be a nonnegative integer. Then we have 
\begin{displaymath}
\kappa_{0,\lambda}(Y)=0,\quad \kappa_{n,\lambda}(Y)=\sum_{k=1}^{n}(-1)^{k-1}\lambda^{k-1}\langle 1\rangle_{k,1/\lambda}{n\brace k}_{Y,\lambda},\quad(n \ge 1).
\end{displaymath}
\end{theorem}
We consider the {\it{probabilistic degenerate Stirling numbers of the first kind associated with $Y$}} defined by
\begin{equation}
\frac{1}{k!}\Big(\log_{-\lambda}\big(E\big[e_{\lambda}^{Y}(t)\big]\big)\Big)^{k}=\sum_{n=k}^{\infty}(-1)^{n-k}S_{1,\lambda}^{Y}(n,k)\frac{t^{n}}{n!},\quad (k\ge 0). \label{21}	
\end{equation}

By using the binomial expansion and \eqref{12}, we have 
\begin{align}
\Big(E\big[e_{\lambda}^{Y}(t)\big]\Big)^{x}&=\Big(E\big[e_{\lambda}^{Y}(t)\big]-1+1\Big)^{x}=\sum_{k=0}^{\infty}(x)_{k}\frac{1}{k!}\Big(E\big[e_{\lambda}^{Y}(t)\big]-1\Big)^{k}\label{22}\\
&=\sum_{k=0}^{\infty}(x)_{k}\sum_{n=k}^{\infty}{n \brace k}_{Y,\lambda}\frac{t^{n}}{n!}=\sum_{n=0}^{\infty}\sum_{k=0}^{n}{n\brace k}_{Y,\lambda}(x)_{k}\frac{t^{n}}{n!}.\nonumber	
\end{align}
On the other hand, by \eqref{6}, we get 
\begin{align}
\Big(E\big[e_{\lambda}^{Y}(t)\big]\Big)^{x}&=e_{-\lambda}^{x}\Big(\log_{-\lambda}\Big(E\big[e_{\lambda}^{Y}(t)\big]\Big)\Big) \label{23} \\
&=\sum_{k=0}^{\infty}(x)_{k,-\lambda}\frac{1}{k!}\Big(\log_{-\lambda}\Big(E\big[e_{\lambda}^{Y}(t)\big]\Big)\Big)^{k} \nonumber \\
&=\sum_{k=0}^{\infty}\langle x\rangle_{k,\lambda}\sum_{n=k}^{\infty}(-1)^{n-k}S_{1,\lambda}^{Y}(n,k)\frac{t^{n}}{n!}\nonumber \\
&=\sum_{n=0}^{\infty}\sum_{k=0}^{n}\langle x\rangle_{k,\lambda}(-1)^{n-k}S_{1,\lambda}^{Y}(n,k)\frac{t^{n}}{n!}.\nonumber	
\end{align}
Therefore, by \eqref{22} and \eqref{23}, we obtain the following theorem. 
\begin{theorem}
For $n\ge 0$, we have 
\begin{displaymath}
\Big(E\big[e_{\lambda}^{Y}(t)\big]\Big)^{x}=\sum_{n=0}^{\infty}\sum_{k=0}^{n}{n\brace k}_{Y,\lambda}(x)_{k}\frac{t^{n}}{n!}=\sum_{n=0}^{\infty}\sum_{k=0}^{n}(-1)^{n-k}\langle x\rangle_{k,\lambda}S_{1,\lambda}^{Y}(n,k)\frac{t^{n}}{n!}. 
\end{displaymath}
\end{theorem}
Applying Theorem 2.2 with $x=m$, we obtain the following corollary. 
\begin{corollary}
For $n\ge 0$, $m\in\mathbb{N}$, we have 
\begin{displaymath}
E\big[(S_{m})_{n,\lambda}\big]=\sum_{k=0}^{n}{n\brace k}_{Y,\lambda}(m)_{k}=\sum_{k=0}^{n}S_{1,\lambda}^{Y}(n,k)(-1)^{n-k}\langle m\rangle_{k,\lambda}.
\end{displaymath}
\end{corollary}
By Theorem 2.2 and \eqref{5}, we get 
\begin{align}
\sum_{k=0}^{n}{n\brace k}_{Y,\lambda}(x)_{k}&=\sum_{k=0}^{n}(-1)^{n-k}(x)_{k,-\lambda}S_{1,\lambda}^{Y}(n,k) \label{24} \\
&=\sum_{k=0}^{n}(-1)^{n-k}S_{1,\lambda}^{Y}(n,k)\sum_{l=0}^{k}{k \brace l}_{-\lambda}(x)_{l} \nonumber \\
&=\sum_{l=0}^{n}\sum_{k=l}^{n}(-1)^{n-k}S_{1,\lambda}^{Y}(n,k){k \brace l}_{-\lambda}(x)_{l}.\nonumber
\end{align}
From \eqref{18} and \eqref{24}, we note that 
\begin{align}
\sum_{k=0}^{n}(-1)^{n-k}S_{1,\lambda}^{Y}(n,k)(x)_{k,-\lambda}&=\sum_{k=0}^{n}{n \brace k}_{Y,\lambda}(-1)^{k}\langle -x\rangle_{k}\label{25} \\
&=\sum_{k=0}^{n}{n\brace k}_{Y,\lambda}(-1)^{k}\sum_{l=0}^{k}{k \brack l}_{-\lambda}\langle -x\rangle_{l,-\lambda}\nonumber\\
&=\sum_{k=0}^{n}{n \brace k}_{Y,\lambda}\sum_{l=0}^{k}(-1)^{k-l}{k \brack l}_{-\lambda}(x)_{l,-\lambda}\nonumber \\
&=\sum_{l=0}^{n}\sum_{k=l}^{n}{n\brace k}_{Y,\lambda}(-1)^{k-l}{k \brack l}_{-\lambda}(x)_{l,-\lambda}\nonumber \\
&=\sum_{l=0}^{n}\sum_{k=l}^{n}{n \brace k}_{Y,\lambda}S_{1,-\lambda}(k,l)(x)_{l,-\lambda}.\nonumber
\end{align}
Therefore, by \eqref{24} and \eqref{25}, we obtain the following theorem. 
\begin{theorem}
For $n\ge l\ge 0$, we have 
\begin{displaymath}
{n \brace l}_{Y,\lambda}=\sum_{k=l}^{n}(-1)^{n-k}S_{1,\lambda}^{Y}(n,k){k \brace l}_{-\lambda},
\end{displaymath}
and 
\begin{displaymath}
(-1)^{n-l}S_{1,\lambda}^{Y}(n,l)=\sum_{k=l}^{n}{n\brace k}_{Y,\lambda}S_{1,-\lambda}(k,l).
\end{displaymath}
\end{theorem}
We note from \eqref{24} and \eqref{17} that 
\begin{align}
\sum_{k=0}^{n}{n\brace k}_{Y,\lambda}(x)_{k}&=\sum_{k=0}^{n}(x)_{k,-\lambda}(-1)^{n-k}S_{1,\lambda}^{Y}(n,k)\label{26}\\
&=\sum_{k=0}^{n}\lambda^{k}(-1)^{n-k}S_{1,\lambda}^{Y}(n,k)\bigg\langle\frac{x}{\lambda}\bigg\rangle_{k}\nonumber\\
&=\sum_{k=0}^{n}\lambda^{k}(-1)^{n-k}S_{1,\lambda}^{Y}(n,k)\sum_{l=0}^{k}L(k,l)\bigg(\frac{x}{\lambda}\bigg)_{l}\nonumber \\
&=\sum_{l=0}^{n}\sum_{k=l}^{n}\lambda^{k-l}(-1)^{n-k}S_{1,\lambda}^{Y}(n,k)L(k,l)(x)_{l,\lambda}\nonumber \\
&=\sum_{l=0}^{n}\sum_{k=l}^{n}\lambda^{k-l}(-1)^{n-k}S_{1,\lambda}^{Y}(n,k)L(k,l)\sum_{j=0}^{l}{l \brace j}_{\lambda}(x)_{j} \nonumber\\
&=\sum_{j=0}^{n}\sum_{l=j}^{n}\sum_{k=l}^{n}\lambda^{k-l}(-1)^{n-k}S_{1,\lambda}^{Y}(n,k)L(k,l){l \brace j}_{\lambda}(x)_{j}.\nonumber
\end{align}
Therefore, by comparing the coefficients on both sides of \eqref{26}, we obtain the following theorem. 
\begin{theorem}
For $n\ge j\ge 0$, we have 
\begin{equation*}
{n \brace j}_{Y,\lambda}=\sum_{l=j}^{n}\sum_{k=l}^{n}\lambda^{k-l}(-1)^{n-k}S_{1,\lambda}^{Y}(n,k)L(k,l){l \brace j}_{\lambda}.
\end{equation*}
\end{theorem}
Taking $k=1$ in \eqref{12}, we have 
\begin{equation}
\sum_{n=1}^{\infty}{n \brace 1}_{Y,\lambda}\frac{t^{n}}{n!}= E\big[e_{\lambda}^{Y}(t)\big]-1=\sum_{n=1}^{\infty}E\big[(Y)_{n,\lambda}\big]\frac{t^{n}}{n!}.\label{27}
\end{equation}
Thus, by \eqref{27}, we get 
\begin{equation}
{n\brace 1}_{Y,\lambda}=E\big[(Y)_{n,\lambda}\big],\quad (n\ge 1). \label{28}	
\end{equation}
From Theorem 2.5 with $j=1$ and using \eqref{28} with $Y=1$, we have 
\begin{align}
{n \brace 1}_{Y,\lambda}&=\sum_{l=1}^{n}\sum_{k=l}^{n}\lambda^{k-l}(-1)^{k-l}S_{1,\lambda}^{Y}(n,k)L(k,l){l \brace 1}_{\lambda}\label{29}\\
&=\sum_{l=1}^{n}\sum_{k=l}^{n}\lambda^{k-l}(-1)^{k-l}(1)_{l,\lambda}S_{1,\lambda}^{Y}(n,k)L(k,l). \nonumber
\end{align}
Therefore, by \eqref{28} and \eqref{29}, we obtain the following theorem. 
\begin{theorem}
For $n\in\mathbb{N}$, we have 
\begin{displaymath}
E[(Y)_{n,\lambda}]=\sum_{l=1}^{n}\sum_{k=l}^{n}\lambda^{k-l}(-1)^{k-l}(1)_{l,\lambda}S_{1,\lambda}^{Y}(n,k)L(k,l).
\end{displaymath}
\end{theorem}
Now, we observe from \eqref{19} and \eqref{21} that, for $k \ge 1$, we have
\begin{align}
&\sum_{n=k}^{\infty}(-1)^{n-k}S_{1,\lambda}^{Y}(n,k)\frac{t^{n}}{n!}=\frac{1}{k!}\Big(\log_{-\lambda}\big(E\big[e_{\lambda}^{Y}(t)\big])\Big)^{k}\label{30}\\
&=\frac{1}{(k-1)!}\Big(\log_{-\lambda}\Big(E\big[e_{\lambda}^{Y}(t)\big]\Big)\Big)^{k-1}\frac{1}{k}\log_{-\lambda}\Big(E\big[e_{\lambda}^{Y}(t)\big]\Big)\nonumber \\
&=\frac{1}{k}\sum_{j=k-1}^{\infty}(-1)^{j-k-1}S_{1,\lambda}^{Y}(j,k-1)\frac{t^{j}}{j!}\sum_{l=1}^{\infty}\kappa_{l,\lambda}(Y)\frac{t^{l}}{l!}\nonumber \\
&=\sum_{n=k}^{\infty}\frac{1}{k}\sum_{j=k-1}^{n-1}\binom{n}{j}(-1)^{j-k-1}S_{1,\lambda}^{Y}(j,k-1)\kappa_{n-j,\lambda}(Y)\frac{t^{n}}{n!}. \nonumber	
\end{align}
Thus, by \eqref{30}, we obtain the following theorem. 
\begin{theorem}
For $n\ge k\ge 1$, we have 
\begin{displaymath}
S_{1,\lambda}^{Y}(n,k)=\frac{1}{k}\sum_{j=k-1}^{n-1}\binom{n}{j}(-1)^{n-j-1}S_{1,\lambda}^{Y}(j,k-1)\kappa_{n-j,\lambda}(Y). 
\end{displaymath}
\end{theorem}
In particular, if $n=k$, then we have 
\begin{align}
S_{1,\lambda}^{Y}(k,k)&=\frac{k}{k}S_{1,\lambda}^{Y}(k-1,k-1)\kappa_{1,\lambda}(Y)=S_{1,\lambda}^{Y}(k-2,k-2)\big(\kappa_{1,\lambda}(Y)\big)^{2} \label{31} \\
&=\big(\kappa_{1,\lambda}(Y)\big)^{3}S_{1,\lambda}^{Y}(k-3,k-3)=\cdots=\big(\kappa_{1,\lambda}(Y)\big)^{k-1}S_{1,\lambda}^{Y}(1,1) \nonumber \\
&=\big(\kappa_{1,\lambda}(Y)\big)^{k}. \nonumber
\end{align}
Therefore, by \eqref{31}, we obtain the following corollary. 
\begin{corollary}
For $k\in\mathbb{N}$, we have 
\begin{displaymath}
S_{1,\lambda}^{Y}(k,k)=\big(\kappa_{1,\lambda}(Y)\big)^{k}.
\end{displaymath}
\end{corollary}
From \eqref{12}, we note that 
\begin{align}
\sum_{n=k}^{\infty}{n \brace k}_{Y,\lambda}\frac{t^{n}}{n!}&=\frac{1}{k!}\Big(E\big[e_{\lambda}^{Y}(t)\big]-1\Big)^{k}\label{32}	\\
&=\frac{1}{(k-1)!}\Big(E\big[e_{\lambda}^{Y}(t)\big]-1\Big)^{k-1}\frac{1}{k}\Big(E\big[e_{\lambda}^{Y}(t)\big]-1\Big) \nonumber \\
&=\sum_{j=k-1}^{\infty}{j \brace k-1}_{Y,\lambda}\frac{t^{j}}{j!}\frac{1}{k}\sum_{l=1}^{\infty}E[(Y)_{l,\lambda}]\frac{t^{l}}{l!} \nonumber \\
&=\sum_{n=k}^{\infty}\frac{1}{k}\sum_{j=k-1}^{n-1}\binom{n}{j}{j \brace k-1}_{Y,\lambda}E[(Y)_{n-j,\lambda}]\frac{t^{n}}{n!},\quad (k \ge 1).\nonumber
\end{align}
Thus, by \eqref{32}, we get 
\begin{equation}
{n \brace k}_{Y,\lambda}=\frac{1}{k}\sum_{j=k-1}^{n-1}\binom{n}{j}{j \brace k-1}_{Y,\lambda}E[(Y)_{n-j,\lambda}],\quad (k \ge 1). \label{33}	
\end{equation}
Taking $n=k$ in \eqref{33}, we have 
\begin{align}
{k \brace k}_{Y,\lambda}&={k-1 \brace k-1}_{Y,\lambda}E[Y]=\Big(E[Y]\Big)^{2}{k-2 \brace k-2}_{Y,\lambda}=\cdots \label{34} \\
&=\Big(E[Y]\Big)^{k-1}{1 \brace 1}_{Y,\lambda}=\Big(E[Y]\Big)^{k}. \nonumber	
\end{align}
Therefore, by \eqref{33} and \eqref{34}, we obtain the following theorem. 
\begin{theorem}
For $n\ge k\ge 1$, we have 
\begin{displaymath}
{n \brace k}_{Y,\lambda}=\frac{1}{k}\sum_{j=k-1}^{n-1}\binom{n}{j}{j\brace k-1}_{Y,\lambda}E[(Y)_{n-j,\lambda}]. 
\end{displaymath}
In particular, we have
\begin{displaymath}
{k \brace k}_{Y,\lambda}=\Big(E[Y]\Big)^{k}. 
\end{displaymath}
\end{theorem}
From \eqref{13} and \eqref{23}, we note that 
\begin{align}
\sum_{n=0}^{\infty}\beta_{n,\lambda}^{Y}(x)\frac{t^{n}}{n!}&=\frac{t}{E\big[e_{\lambda}^{Y}(t)\big]-1}\Big(E\big[e_{\lambda}^{Y}(t)\big]\Big)^{x}\label{35} \\
&=\sum_{l=0}^{\infty}\beta_{l,\lambda}^{Y}\frac{t^{l}}{l!}\sum_{j=0}^{\infty}\sum_{k=0}^{j}\langle x\rangle_{k,\lambda}(-1)^{j-k}S_{1,\lambda}^{Y}(j,k)\frac{t^{j}}{j!} \nonumber \\
&=\sum_{n=0}^{\infty}\sum_{j=0}^{n}\sum_{k=0}^{j}\binom{n}{j}\beta_{n-j,\lambda}^{Y}\langle x\rangle_{k,\lambda}(-1)^{j-k}S_{1,\lambda}^{Y}(j,k)\frac{t^{n}}{n!},\nonumber
\end{align}
where $\beta_{n,\lambda}^{Y}=\beta_{n,\lambda}^{Y}(0)$ are  called the probabilistic degenerate Bernoulli numbers associated with $Y$. By comparing the coefficients on both sides of \eqref{35}, we obtain the following theorem. 
\begin{theorem}
For $n\ge 0$, we have 
\begin{displaymath}
\beta_{n,\lambda}^{Y}(x)=\sum_{j=0}^{n}\sum_{k=0}^{j}\binom{n}{j}\beta_{n-j,\lambda}^{Y}\langle x\rangle_{k,\lambda}(-1)^{j-k}S_{1,\lambda}^{Y}(j,k). 
\end{displaymath}
\end{theorem}
By using \eqref{15} and \eqref{23} and proceeding just as in \eqref{35}, we obtain the next theorem. Here we note that 
$\mathcal{E}_{n,\lambda}^{Y}=\mathcal{E}_{n,\lambda}^{Y}(0)$ are the probabilistic degenerate Euler numbers associated with $Y$. \par  
\begin{theorem}
For $n \ge 0$, we have 
\begin{displaymath}
\mathcal{E}_{n,\lambda}^{Y}(x)=\sum_{j=0}^{n}\sum_{k=0}^{j}\binom{n}{j}\mathcal{E}_{n-j,\lambda}^{Y}\langle x\rangle_{k,\lambda}(-1)^{j-k}S_{1,\lambda}^{Y}(j,k). 
\end{displaymath}
\end{theorem}
By \eqref{15}, we get 
\begin{align}
\Big(E\big[e_{\lambda}^{Y}(t)\big]\Big)^{x}&=\frac{1}{2}\sum_{l=0}^{\infty}\mathcal{E}_{l,\lambda}^{Y}(x)\frac{t^{l}}{l!}\Big(E\big[e_{\lambda}^{Y}(t)\big]+1\Big) \label{36} \\
&=\sum_{n=0}^{\infty}\frac{1}{2}\bigg(\sum_{l=0}^{n}\binom{n}{l}\mathcal{E}_{l,\lambda}^{Y}(x)E\big[(Y)_{n-l,\lambda}\big]+\mathcal{E}_{n,\lambda}^{Y}(x)\bigg)\frac{t^{n}}{n!}.\nonumber	
\end{align}
By using \eqref{7} and \eqref{12}, we note that 
\begin{align}
\Big(E\big[e_{\lambda}^{Y}(t)\big]\Big)^{x}&=e_{\lambda}^{x}\Big(\log_{\lambda}\big(E\big[e_{\lambda}^{Y}(t)\big]-1+1\big)\Big) \label{37} \\
&=\sum_{k=0}^{\infty}(x)_{k,\lambda}\frac{1}{k!}\Big(\log_{\lambda}\big(E\big[e_{\lambda}^{Y}(t)\big]-1+1\big)\Big)^{k} \nonumber \\
&=\sum_{k=0}^{\infty}(x)_{k,\lambda}\sum_{m=k}^{\infty}S_{1,\lambda}(m,k)\frac{1}{m!}\Big(E\big[e_{\lambda}^{Y}(t)\big]-1\Big)^{n} \nonumber \\
&=\sum_{m=0}^{\infty}\sum_{k=0}^{m}(x)_{k,\lambda}S_{1,\lambda}(m,k)\sum_{n=m}^{\infty}{n \brace m}_{Y,\lambda}\frac{t^{n}}{n!}\nonumber\\
&=\sum_{n=0}^{\infty}\sum_{m=0}^{n}\sum_{k=0}^{m}(x)_{k,\lambda}S_{1,\lambda}(m,k){n \brace m}_{Y,\lambda}\frac{t^{n}}{n!} \nonumber \\
&=\sum_{n=0}^{\infty}\sum_{k=0}^{n}(x)_{k,\lambda}\sum_{m=k}^{n}S_{1,\lambda}(n,k){n \brace m}_{Y,\lambda}\frac{t^{n}}{n!}. \nonumber
\end{align}
Therefore, by \eqref{36}, \eqref{37} and Theorem 2.2, we obtain the following theorem. 
\begin{theorem}
For $n\ge 0$, we have 
\begin{align*}
\frac{1}{2}\bigg(\sum_{l=0}^{n}\binom{n}{l}\mathcal{E}_{l,\lambda}^{Y}(x)E\big[(Y)_{n-l,\lambda}\big]+\mathcal{E}_{n,\lambda}^{Y}(x)\bigg)&=\sum_{k=0}^{n}\langle x\rangle_{k,\lambda}(-1)^{n-k}S_{1,\lambda}^{Y}(n,k)\\
&=\sum_{k=0}^{n}(x)_{k,\lambda}\sum_{m=k}^{n}S_{1,\lambda}(m,k){n \brace m}_{Y,\lambda}.
\end{align*}
\end{theorem}

\section{Applications to normal and gamma random variables}
Let $Y\sim N(0,1)$. Then we have 
\begin{align}
E\big[e_{\lambda}^{Y}(t)\big]&=\frac{1}{\sqrt{2\pi}}\int_{-\infty}^{\infty}e^{-\frac{y^{2}}{2}}e_{\lambda}^{y}(t)dy \label{38} \\
&=\frac{1}{\sqrt{2\pi}}\int_{-\infty}^{\infty} e^{-\frac{1}{2}\big(y^{2}-2\frac{y}{\lambda}\log(1+\lambda t)+\big(\frac{\log(1+\lambda t)}{\lambda}\big)^{2}\big)+\frac{1}{2}\big(\frac{\log(1+\lambda t)}{\lambda}\big)^{2}}dy	\nonumber \\
&=\frac{1}{\sqrt{2\pi}}e^{\frac{1}{2}(\log (e_{\lambda}(t))^{2}}\int_{-\infty}^{\infty}e^{-\frac{1}{2}\big(y-\log e_{\lambda}(t)\big)^{2}}dy \nonumber \\
&=e^{\frac{1}{2}(\log e_{\lambda}(t))^{2}}.\nonumber
\end{align}
Thus, by \eqref{38}, we get 
\begin{align}
\sum_{n=0}^{\infty}E\big[(Y)_{n,\lambda}\big]\frac{t^{n}}{n!}&=\sum_{k=0}^{\infty}\frac{(2k)!}{2^{k}k!}\frac{1}{(2k)!}\big(\log(e_{\lambda}(t)-1+1)\big)^{2k} \label{39} \\
&=\sum_{k=0}^{\infty}\frac{(2k)!}{2^{k}k!}\sum_{m=2k}^{\infty}S_{1}(m,2k)\frac{1}{m!}\big(e_{\lambda}(t)-1\big)^{m} \nonumber \\
&=\sum_{m=0}^{\infty}\sum_{k=0}^{[\frac{m}{2}]}\frac{(2k)!}{2^{k}k!}S_{1}(m,2k)\sum_{n=m}^{\infty}{n \brace m}_{\lambda}\frac{t^{n}}{n!} \nonumber \\
&=\sum_{n=0}^{\infty}\sum_{m=0}^{n}\sum_{k=0}^{[\frac{m}{2}]}\frac{(2k)!}{2^{k}k!}S_{1}(m,2k){n \brace m}_{\lambda}\frac{t^{n}}{n!}.\nonumber
\end{align}
Thus, by \eqref{39}, we obtain the following theorem.
\begin{theorem}
Let $Y\sim N(0,1)$. For $n\ge 0$, we have 
\begin{equation*}
E[(Y)_{n,\lambda}]=\sum_{m=0}^{n}\sum_{k=0}^{[\frac{m}{2}]}\frac{(2k)!}{2^{k}k!}S_{1}(m,2k){n \brace m}_{\lambda}=\sum_{k=0}^{[\frac{n}{2}]}\sum_{m=2k}^{n}\frac{(2k)!}{2^{k}k!}S_{1}(m,2k){n \brace m}_{\lambda}.
\end{equation*}
In particular, we have
\begin{displaymath}
E\big[(Y)_{2n,\lambda}\big]=\sum_{k=0}^{n}\sum_{m=2k}^{2n}\frac{(2k)!}{2^{k}k!}S_{1}(m,2k){2n \brace m}_{\lambda}. 
\end{displaymath}
\end{theorem}
Note that 
\begin{align}
E[Y^{2n}]&=\lim_{\lambda\rightarrow 0}E[(Y)_{2n,\lambda}]=\sum_{k=0}^{n}\sum_{m=2k}^{2n}\frac{(2k)!}{2^{k}k!}S_{1}(m,2k){2n \brace m} \label{40} \\
&=\sum_{k=0}^{n}\frac{(2k)!}{2^{k}k!}\sum_{m=2k}^{2n}{n \brace m}S_{1}(m,2k)=\frac{(2n)!}{2^{n}n!}.\nonumber
\end{align}
Let $Y\sim N(0,1)$. Then, from \eqref{38}, we have 
\begin{align}
&\sum_{n=k}^{\infty}{n \brace k}_{Y,\lambda}\frac{t^{n}}{n!}=\frac{1}{k!}\Big(E\big[e_{\lambda}^{Y}(t)\big]-1\Big)^{k}=\frac{1}{k!}\Big(e^{\frac{1}{2}(\log e_{\lambda}(t))^{2}}-1\Big)^{k} \label{41} \\
&=\frac{1}{k!}\sum_{l=0}^{k}\binom{k}{l}(-1)^{k-l}\sum_{j=0}^{\infty}\frac{l^{j}(2j)!}{2^{j}j!}\frac{1}{(2j)!}\Big(\log\big(e_{\lambda}(t)-1+1\big)\Big)^{2j}\nonumber \\
&=\frac{1}{k!}\sum_{l=0}^{k}\binom{k}{l}(-1)^{k-l}\sum_{j=0}^{\infty}\frac{l^{j}(2j)!}{2^{j}j!}\sum_{m=2j}^{\infty}S_{1}(m,2j)\frac{1}{m!}\big(e_{\lambda}(t)-1\big)^{m} \nonumber \\
&=\frac{1}{k!}\sum_{m=0}^{\infty}\sum_{j=0}^{[\frac{m}{2}]}\sum_{l=0}^{k}\binom{k}{l}(-1)^{k-l}\frac{l^{j}(2j)!}{2^{j}j!}S_{1}(m,2j)\sum_{n=m}^{\infty}{n \brace m}_{\lambda}\frac{t^{n}}{n!} \nonumber \\
&=\frac{1}{k!}\sum_{n=0}^{\infty}\sum_{m=0}^{n}\sum_{j=0}^{[\frac{m}{2}]}\sum_{l=0}^{k}\binom{k}{l}(-1)^{k-l}\frac{l^{j}(2j)!}{2^{j}j!}S_{1}(m,2j){n\brace m}_{\lambda}\frac{t^{n}}{n!}\nonumber \\
&=\frac{1}{k!}\sum_{n=0}^{\infty}\sum_{j=0}^{[\frac{n}{2}]}\sum_{m=2j}^{n}
\sum_{l=0}^{k}\binom{k}{l}(-1)^{k-l}\frac{l^{j}(2j)!}{2^{j}j!}S_{1}(m,2j){n\brace m}_{\lambda}\frac{t^{n}}{n!}\nonumber \\
&=\sum_{n=0}^{\infty}\bigg(\frac{1}{k!}\sum_{l=0}^{k}\binom{k}{l}(-1)^{k-l}\sum_{j=0}^{[\frac{n}{2}]}\frac{l^{j}(2j)!}{2^{j}j!}\sum_{m=2j}^{n}S_{1}(m,2j){n\brace m}_{\lambda}\bigg)\frac{t^{n}}{n!}.\nonumber
\end{align}
Therefore, by \eqref{41}, we obtain the following theorem. 
\begin{theorem}
Let $Y\sim N(0,1)$. For $n\ge k\ge 0$, we have 
\begin{displaymath}
{n \brace k}_{Y,\lambda}=\frac{1}{k!}\sum_{l=0}^{k}\binom{k}{l}(-1)^{k-l}\sum_{j=0}^{[\frac{n}{2}]}\frac{l^{j}(2j)!}{2^{j}j!}\sum_{m=2j}^{n}S_{1}(m,2j){n\brace m}_{\lambda}.
\end{displaymath}
\end{theorem}

Let $Y\sim N(0,1)$. Then, by invoking \eqref{38}, we have 
\begin{align}
\sum_{n=0}^{\infty}\mathcal{E}_{n,\lambda}^{Y}(x)\frac{t^{n}}{n!} &=\frac{2}{E\big[e_{\lambda}^{Y}(t)\big]+1}\Big(E\big[e_{\lambda}^{Y}(t)\big]\Big)^{x}=\frac{2}{e^{\frac{1}{2}(\log e_{\lambda}(t))^{2}}+1}e^{\frac{x}{2}(\log e_{\lambda}(t))^{2}} \label{42} \\
&=\sum_{k=0}^{\infty}E_{k}(x)\frac{(2k)!}{2^{k}k!}\frac{1}{(2k)!}\Big(\log(e_{\lambda}(t)-1+1)\Big)^{2k} \nonumber \\
&=\sum_{k=0}^{\infty}E_{k}(x)\frac{(2k)!}{2^{k}k!}\sum_{j=2k}^{\infty}S_{1}(j,2k)\frac{1}{j!}(e_{\lambda}(t)-1)^{j} \nonumber \\
&=\sum_{j=0}^{\infty}\sum_{k=0}^{[\frac{j}{2}]}E_{k}(x)\frac{(2k)!}{2^{k}k!}S_{1}(j,2k)\sum_{n=j}^{\infty}{n\brace j}_{\lambda}\frac{t^{n}}{n!}\nonumber \\
&=\sum_{n=0}^{\infty}\sum_{j=0}^{n}\sum_{k=0}^{[\frac{j}{2}]}E_{k}(x)\frac{(2k)!}{2^{k}k!}S_{1}(j,2k){n\brace k}_{\lambda}\frac{t^{n}}{n!},\nonumber 
\end{align}
where $E_{k}(x)$ are the ordinary Euler polynomials given by
\begin{equation*}
\frac{2}{e^{t}+1}e^{xt}=\sum_{n=0}^{\infty}E_{n}(x)\frac{t^{n}}{n!}. 	
\end{equation*}
Therefore, by \eqref{42}, we obtain the following theorem.
\begin{theorem}
Let $Y\sim N(0,1)$. For $n\ge 0$, we have 
\begin{equation*}
\mathcal{E}_{n,\lambda}^{Y}(x)=\sum_{j=0}^{n}\sum_{k=0}^{[\frac{j}{2}]}E_{k}(x)\binom{2k}{k}\frac{k!}{2^{k}}S_{1}(j,2k){n \brace j}_{\lambda}.
\end{equation*}
\end{theorem}
Let $Y\sim N(0,1)$. Then, by using \eqref{36} and \eqref{38}, we have 
\begin{align}
&\sum_{n=0}^{\infty}\frac{1}{2}\bigg(\sum_{l=0}^{n}\binom{n}{l}\mathcal{E}_{l,\lambda}^{Y}(x)E\big[(Y)_{n-l,\lambda}\big]+\mathcal{E}_{n,\lambda}^{Y}(x)\bigg)\frac{t^{n}}{n!}=\Big(E\big[e_{\lambda}^{Y}(t)\big]\Big)^{x}\label{43} \\
&=e^{\frac{x}{2}(\log(e_{\lambda}(t))^{2}}=\sum_{k=0}^{\infty}\frac{x^{k}(2k)!}{2^{k}k!}\frac{1}{(2k)!}\Big(\log\big(e_{\lambda}(t)-1+1\big)\Big)^{2k}\nonumber \\
&=\sum_{k=0}^{\infty}\frac{(2k)!}{2^{k}k!}x^{k}\sum_{j=2k}^{\infty}S_{1}(j,2k)\frac{1}{j!}\Big(e_{\lambda}(t)-1\Big)^{j}\nonumber \\
&=\sum_{j=0}^{\infty}\sum_{k=0}^{[\frac{j}{2}]}\frac{(2k)!}{2^{k}k!}x^{k}S_{1}(j,2k)\sum_{n=j}^{\infty}{n \brace j}_{\lambda}\frac{t^{n}}{n!}\nonumber \\
&=\sum_{n=0}^{\infty}\sum_{j=0}^{n}\sum_{k=0}^{[\frac{j}{2}]}\frac{(2k)!}{2^{k}k!}S_{1}(j,2k){n \brace j}_{\lambda}x^{k}\frac{t^{n}}{n!}.\nonumber
\end{align}
Therefore, by \eqref{43}, we obtain the following theorem. 
\begin{theorem}
Let $Y\sim N(0,1)$. For $n\ge 0$, we have 
\begin{displaymath}
\frac{1}{2}\bigg(\sum_{l=0}^{n}\binom{n}{l}\mathcal{E}_{l,\lambda}^{Y}(x)E\big[(Y)_{n-l,\lambda}\big]+\mathcal{E}_{n,\lambda}^{Y}(x)\bigg)=\sum_{j=0}^{n}\sum_{k=0}^{[\frac{j}{2}]}\frac{(2k)!}{2^{k}k!}S_{1}(j,2k){n \brace j}_{\lambda}x^{k}.
\end{displaymath}
\end{theorem}
When $Y\sim\Gamma(1,1)$, we have 
\begin{equation}
E\big[e_{\lambda}^{Y}(t)\big]=\int_{0}^{\infty}e_{\lambda}^{y}(t)e^{-y}dy=\frac{1}{1-\frac{1}{\lambda}\log(1+\lambda t)},\quad(t < \frac{1}{\lambda}(e^{\lambda}-1)). \label{44}	
\end{equation}
Thus, by \eqref{21} and \eqref{44}, we get 
\begin{align}
\sum_{n=k}^{\infty}(-1)^{n-k}S_{1,\lambda}^{Y}(n,k)\frac{t^{n}}{n!}&=\frac{1}{k!}\Big(\log_{-\lambda}\big(E\big[e_{\lambda}^{Y}(t)\big]\big)\Big)^{k}\label{45}\\
&=\frac{1}{k!}\bigg(\log_{-\lambda}\frac{1}{1-\frac{1}{\lambda}\log(1+\lambda t)}\bigg)^{k}\nonumber \\
&=\frac{(-1)^{k}}{k!}\Big(\log_{\lambda}\Big(1-\frac{1}{\lambda}\log(1+\lambda t)\Big)\Big)^{k}\nonumber \\
&=(-1)^{k}\sum_{l=k}^{\infty}S_{1,\lambda}(l,k)(-1)^{l}\frac{1}{\lambda^{l}}\frac{1}{l!}\Big(\log(1+\lambda t)\Big)^{l} \nonumber \\
&=(-1)^{k}\sum_{l=k}^{\infty}S_{1,\lambda}(l,k)(-1)^{l}\lambda^{-l}\sum_{n=l}^{\infty}S_{1}(n,l)\frac{\lambda^{n}}{n!}t^{n}\nonumber \\
&=\sum_{n=k}^{\infty}\sum_{l=k}^{n}(-1)^{l-k}S_{1,\lambda}(l,k)S_{1}(n,l)\lambda^{n-l}\frac{t^{n}}{n!}.\nonumber
\end{align}
Therefore, by comparing the coefficients on both sides of \eqref{45}, we obtain the following theorem. 
\begin{theorem}
Let $Y\sim\Gamma(1,1)$. For $n\ge k\ge 0$, we have 
\begin{displaymath}
S_{1,\lambda}^{Y}(n,k)=\sum_{l=k}^{n}(-1)^{n-l}S_{1,\lambda}(l,k)S_{1}(n,l)\lambda^{n-l}. 
\end{displaymath}
In particular, we have
\begin{displaymath}
\lim_{\lambda\rightarrow 0}S_{1,\lambda}^{Y}(n,k)=S_{1}(n,k).
\end{displaymath}
\end{theorem}
By \eqref{13}, we get 
\begin{equation}
\sum_{k=0}^{n}\Big(E\big[e_{\lambda}^{Y}(t)\big]\Big)^{k}=\sum_{m=0}^{\infty}\bigg(\frac{\beta_{m+1,\lambda}^{Y}(n+1)-\beta_{m+1,\lambda}^{Y}}{m+1}\bigg)\frac{t^{m}}{m!}.\label{46}
\end{equation}
Let $Y\sim \Gamma(1,1)$. Then, from \eqref{44}, we have 
\begin{align}
\sum_{k=0}^{n}\Big(E\big[e_{\lambda}^{Y}(t)\big]\Big)^{k}&=\sum_{k=0}^{n}\bigg(\frac{1}{1-\frac{1}{\lambda}\log(1+\lambda t)}\bigg)^{k}	\label{47} \\
&=\sum_{k=0}^{n}\sum_{l=0}^{\infty}\binom{k+l-1}{l}\frac{l!}{\lambda^{l}}\frac{1}{l!}\big(\log(1+\lambda t)\big)^{l} \nonumber \\
&=\sum_{k=0}^{n}\sum_{l=0}^{\infty}\binom{k+l-1}{l}\frac{l!}{\lambda^{l}}\sum_{m=l}^{\infty}S_{1}(m,l)\lambda^{m}\frac{t^{m}}{m!}\nonumber \\
&=\sum_{m=0}^{\infty}\sum_{k=0}^{n}\sum_{l=0}^{m}\binom{k+l-1}{l}l!\lambda^{m-l}S_{1}(m,l)\frac{t^{m}}{m!}.\nonumber
\end{align}
Therefore, by \eqref{46} and \eqref{47}, we obtain the following theorem. 
\begin{theorem}
Let $Y\sim\Gamma(1,1)$. For $n,m\ge 0$, we have 
\begin{displaymath}
\frac{1}{m+1}\Big(\beta_{m+1,\lambda}^{Y}(n+1)-\beta_{m+1,\lambda}^{Y}\Big)=\sum_{k=0}^{n}\sum_{l=0}^{m}\binom{k+l-1}{l}l!\lambda^{m-l}S_{1}(m,l). 
\end{displaymath}
\end{theorem} 

\section{Conclusion} 
Let $Y$ be a random variable whose degenerate moment generating functions exist in some neighborhoods of the origin (see \eqref{10-1}). The probabilistic degenerate Stirling numbers of the second kind associated with $Y$ were studied earlier (see [18,21,22]). Here we studied the probabilistic degenerate Stirling numbers of the first kind, $S_{1,\lambda}^{Y}(n,k)$ (see \eqref{21}), as a degenerate version of the probabilisitic Stirling numbers of the first kind, which were introduced recently by Adell-B\'enyi. We note here that they are constructed from the degenerate cumulant generating function of $Y$ (see \eqref{19}, \eqref{21}). Another thing to note is that $\lim_{\lambda \rightarrow 0} S_{1,\lambda}^{Y}(n,k)=\delta_{n,k}$ (Kronecker's delta) when $Y=1$, while $\lim_{\lambda \rightarrow 0} S_{1,\lambda}^{Y}(n,k)=S_{1}(n,k)$ when $Y \sim \Gamma(1,1)$ (see Theorem 3.5). We investigated some properties, explicit expressions, related identities and recurrence relations for $S_{1,\lambda}^{Y}(n,k)$. In addition, we applied our results to the special cases of normal and gamma random variables. \par
It is one of our future projects to continue to study degenerate versions, $\lambda$-analogues and probabilistic extensions of some special numbers and polynomials and to find their applications to physics, science and engineering as well as to mathematics.

\end{document}